\documentclass[a4paper,11pt]{article}
\usepackage{amsmath,amsthm,amssymb}
\usepackage[bookmarks=false,hyperfootnotes=false]{hyperref}

\numberwithin{equation}{section}

\newcommand{\bh}{\mathcal{B}(H)}

\newcommand{\ot}{\otimes}
\newcommand{\id}{\mathord{\operatorname{id}}}
\newcommand{\Z}{\mathbb{Z}}
\newcommand{\R}{\mathbb{R}}
\newcommand{\N}{\mathbb{N}}
\newcommand{\C}{\mathbb{C}}
\newcommand{\cG}{\mathcal{G}}
\newcommand{\rE}{\operatorname{E}}
\newcommand{\rV}{\operatorname{V}}
\newcommand{\rL}{\mathord{\text{\rm L}}}
\newcommand{\actson}{\curvearrowright}
\newcommand{\recht}{\rightarrow}
\newcommand{\cT}{\mathcal{T}}
\newcommand{\SL}{\operatorname{SL}}
\newcommand{\Stab}{\operatorname{Stab}}
\newcommand{\Aut}{\operatorname{Aut}}

\newcommand{\embed}{\prec}
\newcommand{\rd}{\operatorname{d}}
\newcommand{\cU}{\mathcal{U}}
\newcommand{\eps}{\varepsilon}

\newcommand{\Ptil}{\widetilde{P}}
\newcommand{\m}{\operatorname{m}}
\newcommand{\al}{\alpha}
\newcommand{\be}{\beta}
\newcommand{\M}{\operatorname{M}}
\newcommand{\D}{\operatorname{D}}
\newcommand{\thetatilde}{\widetilde{\theta}}
\newcommand{\om}{\omega}
\newcommand{\vphi}{\varphi}

\newcommand{\bim}[3]{\mathord{\raisebox{-0.4ex}[0ex][0ex]{\scriptsize $#1$}{#2}\hspace{-0.2ex}\raisebox{-0.4ex}[0ex][0ex]{\scriptsize $#3$}}}

\newcommand{\Gr}{\mathcal{G}}

\newcommand{\Mtilde}{\widetilde{M}}
\newcommand{\Ptilde}{\widetilde{P}}

\newcommand{\tree}{\mathcal{T}}
\newcommand{\ET}{\rE(\tree)}
\newcommand{\VT}{\rV(\tree)}

\newcommand{\Gam}{\Gamma}

\newcommand{\Lto}{\rL^2}
\newcommand{\LtoM}{\rL^2(M)}

\newcommand{\HNN}{\operatorname{HNN}}
\theoremstyle{plain}
\newtheorem{theorem}{Theorem}[section]

\newtheorem{lemma}[theorem]{Lemma}
\newtheorem{proposition}[theorem]{Proposition}

\theoremstyle{definition}
\newtheorem{definition}[theorem]{Definition}

\newtheorem{remark}[theorem]{Remark}

\begin{document}
%%%%%%%%%%%%%%%%%%%%%%%%%%%%%%%%%%%%%%%%%%%%%%%%%%%%%%%%%%%%%%%%%%%%%%%%%%%%%%%%%%%%%%%%%%

\begin{center}
{\LARGE\bf HNN extensions and unique group  measure space \vspace{0.5ex}\\ decomposition of II$_1$ factors}

\bigskip

{\sc by Pierre Fima$^{(1,3)}$ and Stefaan Vaes$^{(2,3)}$
\setcounter{footnote}{1}\footnotetext{Supported by ERC Starting Grant VNALG-200749.}
\setcounter{footnote}{2}\footnotetext{Partially supported by ERC Starting Grant VNALG-200749, Research Programme G.0231.07 of the Research Foundation --
Flanders (FWO) and K.U.Leuven BOF research grant OT/08/032.}
\setcounter{footnote}{3}\footnotetext{Department of Mathematics;
    K.U.Leuven; Celestijnenlaan 200B; B--3001 Leuven (Belgium).
    \\ E-mails: pierre.fima@wis.kuleuven.be and stefaan.vaes@wis.kuleuven.be}}

\end{center}

\begin{abstract}\noindent
We prove that for a fairly large family of HNN extensions $\Gamma$, the group measure space II$_1$ factor $\rL^\infty(X) \rtimes \Gamma$ given by an arbitrary free ergodic probability measure preserving action of $\Gamma$, has a unique group measure space Cartan subalgebra up to unitary conjugacy. We deduce from this new examples of W$^*$-superrigid group actions, i.e.\ where the II$_1$ factor $\rL^\infty(X) \rtimes \Gamma$ entirely remembers the group action that it was constructed from.
\end{abstract}

%%%%%%%%%%%%%%%%%%%%%%%%%%%%%%%%%%%%%%%%%%%%%%%%%%%%%%%%%%%%%%%%%%%%%%%%%%%%%%%%%%%%%%%%%%
\section{Introduction and statement of main results}
%%%%%%%%%%%%%%%%%%%%%%%%%%%%%%%%%%%%%%%%%%%%%%%%%%%%%%%%%%%%%%%%%%%%%%%%%%%%%%%%%%%%%%%%%%%%%

A central theme in the theory of von Neumann algebras is the classification of Murray and von Neumann's group measure space II$_1$ factors $\rL^\infty(X) \rtimes \Gamma$ in terms of the initial free ergodic probability measure preserving (p.m.p.) action $\Gamma \actson (X,\mu)$. Typically the II$_1$ factor $\rL^\infty(X) \rtimes \Gamma$ forgets a lot of information about the group action $\Gamma \actson (X,\mu)$. This is best illustrated by Connes' theorem \cite{Co76} implying that all free ergodic p.m.p.\ actions of all infinite amenable groups give rise to the same II$_1$ factor.

Using his groundbreaking deformation/rigidity theory, Popa established several striking rigidity theorems. For particular families of group actions and with different degrees of precision, he manages to recover $\Gamma \actson (X,\mu)$ from the II$_1$ factor $\rL^\infty(X) \rtimes \Gamma$. In particular, in \cite{Po03,Po04} Popa proved that if $\Gamma \actson (X,\mu)$ is any free ergodic p.m.p.\ action of a property (T) group and if $\Lambda \actson (Y,\eta) = [0,1]^\Lambda$ is the Bernoulli action of any group with infinite conjugacy classes (icc), then the isomorphism of $\rL^\infty(X) \rtimes \Gamma$ and $\rL^\infty(Y) \rtimes \Lambda$ implies the isomorphisms of the groups $\Gamma,\Lambda$ and the conjugacy of their actions.

In the recent articles \cite{Pe09,PV09,Io10}, group actions $\Gamma \actson (X,\mu)$ satisfying the most extreme form of rigidity, called W$^*$-superrigidity, were discovered: if $\Lambda \actson (Y,\eta)$ is any free ergodic p.m.p.\ action of any group $\Lambda$ and if $\rL^\infty(X) \rtimes \Gamma$ is isomorphic with $\rL^\infty(Y) \rtimes \Lambda$, then the groups $\Gamma,\Lambda$ are isomorphic and their actions conjugate.

Note that W$^*$-superrigidity for $\Gamma \actson (X,\mu)$ arises as the sum of the following two rigidity phenomena.
\begin{itemize}
\item Uniqueness of the group measure space Cartan subalgebra: if the II$_1$ factor $M = \rL^\infty(X) \rtimes \Gamma$ has another group measure space decomposition $M = \rL^\infty(Y) \rtimes \Lambda$, then the Cartan subalgebras $\rL^\infty(X)$ and $\rL^\infty(Y)$ of $M$ must be unitarily conjugate inside $M$.
\item Orbit equivalence superrigidity: if any other free ergodic p.m.p.\ action $\Lambda \actson (Y,\eta)$ is orbit equivalent with $\Gamma \actson (X,\mu)$, then the groups must be isomorphic and the actions conjugate.
\end{itemize}
Both the uniqueness of (group measure space) Cartan subalgebras and the orbit equivalence superrigidity are extremely hard to establish and are even more rarely known to hold simultaneously.

In \cite{PV09} a class $\cG$ of countable groups $\Gamma$ was found such that all group measure space II$_1$ factors $\rL^\infty(X) \rtimes \Gamma$, for \emph{arbitrary} free ergodic p.m.p.\ actions $\Gamma \actson (X,\mu)$, have a unique group measure space Cartan subalgebra up to unitary conjugacy. This class $\cG$ contains all free products $\Gamma_1 * \Gamma_2$ of an infinite property (T) group $\Gamma_1$ and a non-trivial group $\Gamma_2$, as well as a fairly large family of amalgamated free product groups $\Gamma_1 *_\Sigma \Gamma_2$.

The main aim of this paper is to show the following similar uniqueness of group measure space Cartan subalgebras for crossed products with certain HNN extensions $\Gamma = \HNN(H,\Sigma,\theta)$ where $\Sigma < H$ is an amenable subgroup, $\theta : \Sigma \recht H$ is an injective group homomorphism and $H$ has a certain rigidity property, e.g.\ property (T). Recall that $\HNN(H,\Sigma,\theta)$ is the group generated by a copy of $H$ and an extra generator $t$, called stable letter, with relations $t \sigma t^{-1} = \theta(\sigma)$ for all $\sigma \in \Sigma$.

\begin{theorem}[{See Theorem \ref{UniqueCartanHNN-amplif}}]\label{UniqueCartanHNN}
Let $H$ be a group that contains a non-amenable subgroup with the relative property (T) or that contains two commuting non-amenable subgroups. Let $\Sigma < H$ be an amenable subgroup and $\theta : \Sigma  \recht H$ an injective group homomorphism. Denote by $\Gamma = \HNN(H,\Sigma,\theta)$ the corresponding HNN extension and assume that there exist $g_1,\ldots,g_n \in \Gamma$ such that $\bigcap_{i=1}^n g_i \Sigma g_i^{-1}$ is finite.

For all free ergodic p.m.p.\ actions $\Gam \actson (X,\mu)$, the II$_1$ factor $\rL^\infty(X) \rtimes \Gamma$ has, up to unitary conjugacy, a unique group measure space Cartan subalgebra.
\end{theorem}

Note that the finiteness of $\bigcap_{i=1}^n g_i \Sigma g_i^{-1}$ is automatic when $\Sigma \cap \theta(\Sigma)$ is finite.

To prove Theorem \ref{UniqueCartanHNN} we follow the same strategy as in \cite{PV09}, replacing the length deformation on amalgamated free products and its dilation (see \cite{IPP05}) by an analogous deformation and dilation for HNN extensions. For this we make use of HNN extensions of tracial von Neumann algebras.

Next we observe that those groups that admit a decomposition either as an amalgamated free product over an amenable group or as an HNN extension over an amenable subgroup, are exactly those groups that admit an action on a tree with at least one amenable edge stabilizer. To make a more precise statement, also avoiding trivial amalgamated free product decompositions, we introduce the following terminology. We follow Serre's conventions \cite{Se83} so that a graph $\cG$ is a pair of sets $\rV(\cG),\rE(\cG)$, whose elements are called the vertices, resp.\ edges of $\cG$, equipped with a source map $s : \rE(\cG) \recht \rV(\cG)$ and a range map $r : \rE(\cG) \recht \rV(\cG)$, as well as an involution $\rE(\cG) \recht \rE(\cG)$ called inversion and satisfying the following two properties: $\overline{e} \neq e$ and $s(\overline{e}) = r(e)$ for all $e \in \rE(\cG)$. A tree is a connected graph without cycles. An action on a tree is said to be without inversion if $g \cdot e \neq \overline{e}$ for all $e \in \rE(\cG)$ and all $g \in \Gamma$.

A combination of \cite[Theorem 1.1]{PV09} and Theorem \ref{UniqueCartanHNN} then yields the following.

\begin{theorem}\label{UniqueCartanTrees}
Let $\Gamma$ be a group satisfying the following two properties.
\begin{enumerate}
\item $\Gamma$ contains a non-amenable subgroup with the relative property $(T)$ or $\Gamma$ contains two commuting non-amenable subgroups.

\item $\Gamma$ admits an action $\Gamma \actson \cT$ without inversion on a tree $\cT$ such that there exists a finite subtree with a finite stabilizer and such that there exists an edge $e \in \rE(\cT)$ with the properties that $\Stab e$ is amenable and that the smallest subtrees containing all vertices $\Gamma \cdot s(e)$, resp.\ $\Gamma \cdot r(e)$, are both equal to the whole of $\cT$.
\end{enumerate}

Then, for all free ergodic p.m.p.\ actions $\Gam \actson (X,\mu)$, the II$_1$ factor $\rL^\infty(X) \rtimes \Gamma$ has, up to unitary conjugacy, a unique group measure space Cartan subalgebra.
\end{theorem}

In combination with Popa's cocycle superrigidity theorems \cite[Theorem 0.1]{Po05} and \cite[Theorem 1.1]{Po06a}, we obtain the following new examples of W$^*$-superrigid actions. Recall that if a countable group $\Gamma$ acts on a countable set $I$, the action $\Gamma \actson (X_0,\mu_0)^I$ given by $(g \cdot x)_i = x_{g^{-1} \cdot i}$ is called a generalized Bernoulli action.

All generalized Bernoulli actions are p.m.p. They are ergodic if and only if $\Gamma \cdot i$ is infinite for all $i \in I$. The criterion for essential freeness depends on whether $(X_0,\mu_0)$ has atoms or not: if $(X_0,\mu_0)$ has no atoms, essential freeness is equivalent with every $g \neq e$ acting non-trivially on $I$; if $(X_0,\mu_0)$ has atoms, essential freeness is equivalent with every $g \neq e$ moving infinitely many $i \in I$. We will always implicitly assume that our generalized Bernoulli actions are essentially free.

\begin{theorem}[{See Theorem \ref{thm.Wstarsuperrigid-general}}]\label{thm.Wstarsuperrigid}
Let $H$ be a property (T) group, $\Sigma < H$ an infinite amenable subgroup and $\theta : \Sigma \recht H$ an injective group homomorphism satisfying $\Sigma \cap \theta(\Sigma) = \{e\}$. If $\Gamma \actson I$ such that $\Sigma \cdot i$ is infinite for all $i \in I$, then the generalized Bernoulli action $\Gamma \actson (X_0,\mu_0)^I$ is W$^*$-superrigid.
\end{theorem}

Observe that explicit examples for Theorem \ref{thm.Wstarsuperrigid} can be obtained by considering two different copies of $\Z$ in $\SL_n(\Z)$, $n\geq 3$.

{\bf Acknowledgment.} We are grateful to the referee for pointing us towards \cite{Ue07} where HNN extensions of von Neumann algebras are, up to amplifications, shown to be isomorphic to certain amalgamated free products. As we explain in Remark \ref{rem.ueda}, this does not allow to directly deduce our main results from the analogous results in \cite{PV09}, but this can be used as the starting point of an alternative, slightly more technical proof.

%%%%%%%%%%%%%%%%%%%%%%%%%%%%%%%%%%%%%%%%%%%%%%%%%%%%%%%%%%%%%%%%%%%%%%%%%%%%%%%%%%%
\section{Preliminaries}
%%%%%%%%%%%%%%%%%%%%%%%%%%%%%%%%%%%%%%%%%%%%%%%%%%%%%%%%%%%%%%%%%%%%%%%%%%%%%%%%%%%%%%

%%%%%%%%%%%%
\subsection{Intertwining by bimodules}
%%%%%%%%%%%%%%%

To fix notations, we briefly recall the \emph{intertwining-by-bimodules} technique from \cite[Section 2]{Po03}.

Let $(M,\tau)$ be a tracial von Neumann algebra and assume that $A,B \subset M$ are possibly non-unital von Neumann subalgebras. Denote their respective units by $1_A$ and $1_B$. Then, the following three conditions are equivalent.
\begin{itemize}
\item $1_A \rL^2(M) 1_B$ admits an $A$-$B$-subbimodule that is finitely generated as a right $B$-module.
\item There exist non-zero projections $p \in A$, $q \in B$, a normal unital $*$-homomorphism $\vphi : pAp \recht qBq$ and a non-zero partial isometry $v \in p M q$ satisfying $a v = v \vphi(a)$ for all $a \in pAp$.
\item There is no sequence of unitaries $u_n \in \cU(A)$ satisfying $\|E_B(x u_n y^*)\|_2 \recht 0$ for all $x,y \in 1_B M 1_A$.
\end{itemize}
If one of these equivalent conditions hold, we write $A \embed_M B$. Otherwise, we write $A \not\embed_M B$.

When $M$ is a II$_1$ factor and $A,B \subset M$ are \emph{Cartan subalgebras,} then $A \embed_M B$ if and only if there exists a unitary $u \in \cU(M)$ such that $A = u B u^*$, see \cite[Theorem A.1]{Po01}.

%%%%%%%%%%%%%%%%%%%%%%%%%%%%%%%%%%%%%%%%%%%%%%%%%%%%%%%%%%%%%%%%%%%%%%%%%%%%%%%%%%%%%%
\subsection{Amalgamated free products of tracial von Neumann algebras}\label{subsec.amal}
%%%%%%%%%%%%%%%%%%%%%%%%%%%%%%%%%%%%%%%%%%%%%%%%%%%%%%%%%%%%%%%%%%%%%%%%%%%%%%%%%%%%%%

Let $P_1$, $P_2$ be von Neumann algebras equipped with a faithful normal tracial state $\tau$. Assume that $N$ is a common von Neumann subalgebra of $P_1$ and $P_2$ and that the traces of $P_1$, $P_2$ coincide on $N$. We denote by $P = P_1 *_N P_2$ the amalgamated free product with respect to the unique trace preserving conditional expectations (see \cite{Po91} and \cite{VDN92}).

We record the following result from \cite{IPP05}. Recall first that the quasi-normalizer of a von Neumann subalgebra $Q \subset P$ is defined as the von Neumann algebra generated by the elements $x \in P$ for which there exist $x_1,\ldots,x_n,y_1,\ldots,y_m \in P$ satisfying
$$x Q \subset \sum_{i=1}^n Q x_i \quad\text{and}\quad Q x \subset \sum_{j=1}^m y_j Q \; .$$

\begin{lemma}[Theorem 1.1 in \cite{IPP05}]\label{lemma.IPP}
Let $P = P_1 *_N P_2$ be an amalgamated free product w.r.t.\ trace preserving conditional expectations, as above. Let $p \in P_1$ be a non-zero projection and $Q \subset p P_1 p$ a von Neumann subalgebra of $p P_1 p$ such that $Q \not\embed_{P_1} N$. Then every $Q$-$P_1$-subbimodule $H$ of $p \rL^2(P)$ which is finitely generated as a right $P_1$-module, is contained in $\rL^2(P_1)$. In particular, the quasi-normalizer of $Q$ inside $p P p$ is contained in $p P_1 p$.
\end{lemma}

We also mention the following fact that can be proven by directly applying the third characterization of the intertwining relation $A \embed B$. If $P = P_1 *_N P_2$ is an amalgamated free product and $Q \subset p P_1 p$ is a von Neumann subalgebra satisfying $Q \not\embed_{P_1} N$, then $Q \not\embed_P P_2$.

%%%%%%%%%%%%%%%%%%%%%%%%%%%%%%%%%%%%%%%%%%%%%%%%%%%%%%%%%%%%%%%%%%%%%%%%%%%%%%%%%%%%%%
\subsection{HNN extensions of groups: some notations}\label{HNNgrp}
%%%%%%%%%%%%%%%%%%%%%%%%%%%%%%%%%%%%%%%%%%%%%%%%%%%%%%%%%%%%%%%%%%%%%%%%%%%%%%%%%%%%%%%

Let $H$ be group, $\Sigma<H$ a subgroup and $\theta\,:\,\Sigma\rightarrow H$ an injective group homomorphism. Define the HNN extension $\Gam=\HNN(H,\Sigma,\theta)=\langle H,t\,|\,\theta(\sigma)=t\sigma t^{-1},\,\,\forall \sigma\in\Sigma\rangle$. For $\eps\in\{-1,1\}$ we define
$$\Sigma_{\eps}=\begin{cases} \Sigma &\;\;\text{if}\;\; \eps=1 \; , \\ \theta(\Sigma) &\;\;\text{if}\;\; \eps = -1 \; .\end{cases}$$

We call $g=g_0t^{\eps_1}g_1t^{\eps_2}\cdots t^{\eps_n}g_n\in\Gam$, with $g_i\in H$ and $\eps_i\in\{-1,1\}$, a \textit{reduced expression} if
$g_i \in H - \Sigma_{\eps_i}$ whenever $\eps_i \neq \eps_{i+1}$. By convention, if $n = 0$, the reduced expressions are defined to be the elements $g_0 \in H - \{e\}$. Observe that any element $g\in\Gam-\{e\}$ admits a reduced expression and that the natural number $n$, as well as the sequence $\eps_1,\ldots,\eps_n$ appearing in a reduced expression for $g$, only depend on $g$. We call $n$ the \textit{length} of $g$ and we denote $|g| := n$.

The Bass-Serre tree $\tree$ of the HNN extension $\Gam$ is defined as follows.
$$\VT=\Gam/H\quad\text{and}\quad \ET^+=\Gam/\Sigma,$$
where $\VT$ denotes the set of vertices of $\tree$ and $\ET^+$ denotes the set of positive oriented edges of $\tree$. The source map $s$ and the range map $r$ are defined by
$$s(g\Sigma)=gH\quad\text{and}\quad r(g\Sigma)=gt^{-1}H\quad\text{for all}\;\; g\in\Gam \; .$$
We have a natural action of $\Gam$ on $\tree$. Denote by $\rd_{\tree}$ the geodesic distance on $\rV(\cT)$. One checks that
$$\rd_{\tree}(H,g\cdot H)=|g|\quad\text{for all}\;\; g\in\Gam \; .$$

%%%%%%%%%%%%%%%%%%%%%%%%%%%%%%%%%%%%%%%%%%%%%%%%%%%%%%%%%%%%%%%%%%%%%%%%%%%%%%%%
\section{HNN extensions of von Neumann algebras and their length deformation}\label{HNNVN}
%%%%%%%%%%%%%%%%%%%%%%%%%%%%%%%%%%%%%%%%%%%%%%%%%%%%%%%%%%%%%%%%%%%%%%%%%%%%%%%%%

HNN extensions of general von Neumann algebras were introduced in \cite{Ue04}. In this section we follow a different approach to specifically define HNN extensions of finite von Neumann algebras. This will serve as an important technical tool in the proof of Theorem \ref{UniqueCartanTrees}.

Consider first the group case. So, let $H$ be a group, $\Sigma<H$ a subgroup and $\theta\,:\,\Sigma\rightarrow H$ an injective group homomorphism. Associated to this data is a triple $(M,N,\theta)$ where $N=\rL \Sigma$ is viewed as a von Neumann subalgebra of $M=\rL H$ and $\theta$ provides a trace preserving embedding from $N$ into $M$.
We shall associate to an arbitrary triple $(M,N,\theta)$ consisting of a tracial von Neumann algebra $M$, a von Neumann subalgebra $N \subset M$ and a trace preserving embedding $\theta : N \recht M$, a new tracial von Neumann algebra $\HNN(M,N,\theta)$. By construction we will have $\rL(\HNN(H,\Sigma,\theta)) = \HNN(\rL H,\rL \Sigma,\theta)$.

Fix such a triple $(M,N,\theta)$.

\subsection{The $\Lto$-space}

For $\eps\in\{-1,1\}$ define
$$N_{\eps}=\begin{cases} N &\;\;\text{if}\;\; \eps = 1 \, ,\\ \theta(N) &\;\;\text{if}\;\; \eps=-1 \; . \end{cases}$$
We define $\theta^\eps : N_\eps \recht N_{-\eps} \subset M$ in the obvious way.

For $n\geq 1$ and $\eps_1,\ldots,\eps_n\in\{-1,1\}$ define the $M$-$M$-bimodule
$$H_{\eps_1,\ldots,\eps_n}=K_0\underset{N}{\otimes}\cdots\underset{N}{\otimes}K_n,$$
where $K_0=K_n=\LtoM$ and, for $1\leq i\leq n-1$,
$$K_i=\begin{cases}
\LtoM &\;\;\text{if}\;\; \eps_i=\eps_{i+1} \; ,\\
\LtoM\ominus \Lto(N_{\eps_i}) &\;\;\text{if}\;\; \eps_i\neq\eps_{i+1}\; .\end{cases}$$
We view $K_0=\LtoM$  as an $M$-$N$-bimodule, where the left $M$-action is the obvious one and the right $N$-action is given by
$\xi\cdot x=\xi x$ if $\eps_1=-1$ and  $\xi\cdot x=\xi \theta(x)$ if $\eps_1=1$. Similarly, we view $K_n=\LtoM$  as an $N$-$M$-bimodule, where the right $M$-action is the obvious one and the left $N$-action is given by
$x\cdot\xi=x\xi$ if $\eps_n=1$ and  $x\cdot\xi=\theta(x)\xi$ if $\eps_n=-1$. Finally, for $1\leq i\leq n-1$, we view $K_i$ as an $N$-$N$-bimodule in the following way.
\begin{itemize}
\item The left $N$-action is given by $x\cdot\xi=\begin{cases}
x\xi &\;\;\text{if}\;\; \eps_i=1 \; ,\\
\theta(x)\xi &\;\;\text{if}\;\;\eps_i=-1 \; .\end{cases}$
\item The right $N$-action is given by $\xi\cdot x= \begin{cases}
\xi\theta(x) &\;\;\text{if}\;\; \eps_{i+1}=1 \; ,\\
\xi x &\;\; \text{if}\;\; \eps_{i+1}=-1 \; .\end{cases}$
\end{itemize}
Define the $M$-$M$-bimodule
$$H=\LtoM\oplus\bigoplus_{n\geq 1,\,\,\eps_{1},\ldots,\eps_n\in\{-1,1\}}H_{\eps_1,\ldots,\eps_n}.$$
We always view $M\subset\bh$ via the left module action and we denote by $\rho\,:\,M\rightarrow\bh$ the unital normal $*$-anti-homomorphism given by the right action of $M$ on $H$. Every element $x \in M$ can be viewed as a vector in $\rL^2(M)$ that we denote by $\hat{x}$.

Let $\eps\in\{-1,1\}$. We define a unitary $u^{\eps}\in\bh$ in the following way.
\begin{itemize}
\item If $\xi\in\LtoM$ we define $u^{\eps}\xi=\hat{1}\ot\xi\in H_{\eps}$.
\item If $\xi\in H_{\eps_1,\ldots,\eps_n}$ with $n\geq 1$ and $\eps_1=\eps$ we define $u^{\eps}\xi=\hat{1}\ot\xi\in H_{\eps,\eps_1,\ldots,\eps_n}$.
\item If $\xi=\hat{x}\ot\xi_0\in H_{\eps_1,\ldots,\eps_n}$ with $n\geq 1$, $\eps_1\neq\eps$ and $x\in M$, $\xi_0\in H_{\eps_2,\ldots,\eps_n}$ we define
$$u^{\eps}(\hat{x}\ot\xi_0)=\begin{cases}
\hat{1}\ot\hat{x}\ot\xi_0 \in H_{\eps,\eps_1,\ldots,\eps_n} &\;\; \text{if}\;\; x\in M\ominus N_{\eps} \; ,\\
\theta^\eps(x)\xi_0 \in H_{\eps_2,\ldots,\eps_n} & \;\;\text{if}\;\; x\in N_{\eps}\; .\end{cases}$$
\end{itemize}
It is easy to check that $u^{\eps}$ extends to a unitary on $H$ such that $(u^{\eps})^*=u^{-\eps}$, justifying the superscript notation. We rather write $u$ instead of $u^1$ and one checks easily that
$$uxu^*=\theta(x)\quad\text{for all}\quad x\in N.$$
Moreover, $u$ commutes with the right $M$-module action on $H$.

In an entirely similar way, we define the right version $v^\eps$ of $u^\eps$.

\begin{itemize}
\item If $\xi\in\LtoM$ we define $v^{\eps}\xi=\xi\ot\hat{1}\in H_{\eps}$.
\item If $\xi\in H_{\eps_1,\ldots,\eps_n}$ with $n\geq 1$ and $\eps_n=\eps$ we define $v^{\eps}\xi=\xi\ot\hat{1}\in H_{\eps_1,\ldots,\eps_n,\eps}$.
\item If $\xi=\xi_0\ot\hat{x}\in H_{\eps_1,\ldots,\eps_n}$ with $n\geq 1$, $\eps_n\neq\eps$ and $x\in M$, $\xi_0\in H_{\eps_1,\ldots,\eps_{n-1}}$ we define
$$v^{\eps}(\xi_0\ot\hat{x})=\begin{cases}
\xi_0\ot\hat{x}\ot \hat{1} \in H_{\eps_1,\ldots,\eps_n,\eps} &\;\;\text{if}\;\; x\in M\ominus N_{-\eps} \; ,\\
\xi_0\theta^{-\eps}(x) \in H_{\eps_1,\ldots,\eps_{n-1}} &\;\;\text{if}\;\; x\in N_{-\eps}\; .\end{cases}$$
\end{itemize}
As above, $v^{\eps}$ extends to a unitary on $H$ such that $(v^{\eps})^*=v^{-\eps}$. Again we write $v$ instead of $v^1$. One checks that $v$ commutes with the left $M$-module action on $H$, as well as with the unitary $u$.

\subsection{The HNN extension}

\begin{definition}
The HNN extension $\HNN(M,N,\theta)$ is defined as the von Neumann subalgebra of $\bh$ generated by $M$ and $u$:
$$\HNN(M,N,\theta):=\langle M,u\rangle\subset\bh.$$
\end{definition}

Let $P=\HNN(M,N,\theta)$. An element $x\in P$ of the form $x=x_0u^{\eps_1}x_1\cdots u^{\eps_n}x_n$ with $x_i\in M$ and $\eps_i\in\{-1,1\}$ will be called \textit{reduced} if $x_i \in M \ominus N_{\eps_i}$ whenever $\eps_i \neq \eps_{i+1}$. By convention, in the case where $n=0$, the reduced elements are the ones of the form $x = x_0$ with $x_0 \in M \ominus \C 1$.

Let $\Omega=\hat{1}\in\LtoM\subset H$. Let $x=x_0u^{\eps_1}x_1\cdots u^{\eps_n}x_n\in P$ be a reduced element. Observe that
$$x\Omega=\hat{x}_0\ot\cdots\ot\hat{x}_n\in H_{\eps_1,\ldots,\eps_n}.$$
It follows that the integer $n$, as well as the sequence $\eps_1,\ldots,\eps_n$, only depend on $x$. The integer $n$ is called the \textit{length} of the reduced element $x$.

Let $P_{{\rm red}}$ be the vector subspace of $P$ spanned by the reduced elements. By the relation $\theta(x)=uxu^*$ for $x\in N$, it is easy to check that $P_{{\rm red}}$ is a $*$-subalgebra of $P$. Moreover, by definition of the HNN extension, $P_{{\rm red}}$ is weakly dense in $P$.

We can also consider the right version of the HNN extension
$$\HNN'(M,N,\theta):=\langle \rho(M),v\rangle\subset\bh.$$
Observe that $\HNN'(M,N,\theta)\subset P'$. It will become clear immediately that this inclusion actually is an equality.

As before, we have the notion a of reduced element $x'=\rho(x_0)v^{\eps_1}\rho(x_1)\cdots v^{\eps_n}\rho(x_n)$ in $\HNN'(M,N,\theta)$ and the vector space spanned by such elements is a weakly dense $*$-subalgebra of $\HNN'(M,N,\theta)$.

\subsection{The HNN trace}

Let $\Omega=\hat{1}\in\LtoM\subset H$. Define the normal state on $P=\HNN(M,N,\theta)$ given by $\tau(x)=\langle \Omega,x \Omega\rangle$. Observe that, whenever $x=x_0u^{\eps_1}x_1\cdots u^{\eps_n}x_n\in P$ is reduced, also $x'=\rho(x_n)v^{\eps_n}\cdots v^{\eps_1}\rho(x_0)$ is a reduced element in $\HNN'(M,N,\theta)$ and one has
$$x\Omega=\hat{x}_0\ot\cdots\ot\hat{x}_n=x'\Omega.$$
It follows that $\Omega$ is a cyclic vector for both $\HNN(M,N,\theta)$ and $\HNN'(M,N,\theta)$. Hence, $\tau$ is a faithful state on $P$ and $(H,\Omega)$ is its GNS construction. Denote by $J$ the modular conjugation on $\rL^2(M)$ given by $J \hat{x} = \widehat{x^*}$. Define the anti-unitary operator $J_\tau$ on $H$ such that ${J_\tau}_{|\rL^2(M)} = J$ and ${J_\tau}_{| H_{\eps_1,\ldots,\eps_n}} = J_{\eps_1,\ldots,\eps_n}$,
where $J_{\eps_1,\ldots,\eps_n}\,:\,H_{\eps_1,\ldots,\eps_n}\rightarrow H_{\eps_n,\ldots,\eps_1}$ is given by the formula $\xi_0\ot\xi_1\ot\cdots\ot\xi_n\mapsto J\xi_n\ot J\xi_{n-1}\ot\cdots\ot J\xi_0$. One checks straightforwardly that $J_\tau(x \Omega) = x^* \Omega$ for all $x \in P$. Since $J_\tau$ is anti-unitary, it follows that $\tau$ is a trace.

Note that $\tau(x)=0$ whenever $x$ is a reduced element. Also observe that the canonical inclusion $M\subset P$ is trace preserving.
Since clearly $J_\tau \HNN(M,N,\theta) J_\tau = \HNN'(M,N,\theta)$, it follows that $\HNN(M,N,\theta)$ and $\HNN'(M,N,\theta)$ are each other's commutant.

\subsection{The universal property}

We record the following elementary proposition and leave the proof to the reader.

\begin{proposition}\label{prop.universal}
Let $P=\HNN(M,N,\theta)$ be an HNN extension. Assume that $(Q,\tau_Q)$ is any tracial von Neumann algebra, that $\pi\,:\,M\rightarrow Q$ is a trace-preserving embedding and that $w\in Q$ is a unitary satisfying
\begin{itemize}
\item $\pi(\theta(x))=w\pi(x)w^*$ for all $x\in N$,
\item for all reduced $x=x_0u^{\eps_1}\cdots u^{\eps_n}x_n\in P$, we have $\tau_Q(\pi(x_0)w^{\eps_1}\cdots w^{\eps_n}\pi(x_n))=0$.
\end{itemize}
Then there exists a unique trace-preserving $*$-homomorphism $\tilde{\pi}: P \recht Q$ extending $\pi$ and satisfying $\tilde{\pi}(u) = w$.
\end{proposition}

Either using the universal property or checking definitions, one observes that $\HNN(\rL H,\rL \Sigma,\theta) = \rL(\HNN(H,\Sigma,\theta))$ whenever $\Sigma < H$ is a subgroup and $\theta : \Sigma \recht H$ is an injective group homomorphism.

\subsection{The length deformation of an HNN extension} \label{sec.deform}

In \cite{IPP05,PV09} a crucial role is played by the length deformation on an amalgamated free product. Since we need a similar deformation for HNN extensions, we first recall the amalgamated free product case. We call reduced element of $P := P_1 *_N P_2$ any product $x = x_1 \cdots x_n$ where the factors $x_i$ belong alternatingly to $P_1 \ominus N$ and $P_2 \ominus N$. Let $0 < \rho < 1$. In \cite{IPP05}, it is shown that the formula
$$\psi_\rho(x) = \rho^n \, x \quad\text{for all reduced}\;\; x = x_1 \cdots x_n \; ,$$
yields a well defined normal unital completely positive map $\psi_\rho : P \recht P$. If $\rho \recht 1$, then $\psi_\rho \recht \id$ pointwise in $\|\,\cdot\,\|_2$.

One of the main technical ingredients of \cite{IPP05} is the following: whenever  $Q \subset P$ is a von Neumann subalgebra with the relative property (T), there exists $i \in \{1,2\}$ such that $Q \embed_P P_i$. By the relative property (T), we know that $\psi_\rho \recht \id$ in $\|\,\cdot\,\|_2$, uniformly on the unit ball of $P$. In order to deduce from this that actually $Q \embed_P P_i$ for some $i \in \{1,2\}$, a dilation of $\psi_\rho$ is introduced in \cite{IPP05}. One constructs a larger tracial von Neumann algebra $\Ptil \supset P$ together with a continuous one-parameter group of automorphisms $\al_t \in \Aut(\Ptil)$ such that
$$\psi_{\rho_t}(x) = E_P(\al_t(x)) \quad\text{for all}\;\; x \in P \; .$$
Here $\rho_t \recht 1$ when $t \recht 0$. In order to allow for spectral gap rigidity, it is crucial that $\bim{P}{\rL^2(\Ptil \ominus P)}{P}$ is weakly contained in the coarse $P$-$P$-bimodule. This is exactly the case if we amalgamate over an amenable subalgebra $N$.

We now perform similar constructions for the HNN extension $P := \HNN(M,N,\theta)$, leading to well defined normal unital completely positive maps $\m_\rho : P \recht P$ satisfying
\begin{equation}\label{eq.defmrho}
\m_\rho(x) = \rho^n \, x \quad\text{whenever}\;\; x = x_0 u^{\eps_1} x_1 \cdots u^{\eps_n} x_n \;\;\text{is a reduced element in $P$}\; .
\end{equation}

Define $\Mtilde := M *_N (N \ot \rL \Z)$ and denote by $v \in \rL \Z$ the canonical unitary generator. Define the automorphism $\be \in \Aut \Mtilde$ given by $\be(x) = x$ for all $x \in M$ and $\be(v) = v^*$. Let $a \in \rL \Z$ be the unique self-adjoint element with spectrum $[-\pi,\pi]$ such that $v = \exp(i a)$. Note that $\be(a) = - a$ and define $v_s := \exp (is a)$. It follows that $v_s$ is a continuous one-parameter group of unitaries in $\rL \Z$ and that $\be(v_s) = v_{-s}$ for all $s \in \R$. Finally, putting $\rho_s := \sin(\pi s) / (\pi s)$, one checks that $\tau(v_s) = \tau(v_s^*) = \rho_s$ for all $s \in \R$.

We naturally have $M \subset \Mtilde$, so that we can define $\Ptilde := \HNN(\Mtilde,N,\theta)$. Note that $P \subset \Ptilde$ and that the \lq stable unitary letter\rq\ $u \in P$ also serves as stable unitary letter for $\Ptilde$. For all $s \in \R$, the unitary $u v_s \in \Ptilde$ satisfies $u v_s x (u v_s)^* = \theta(x)$ for all $x \in N$. We want to apply the universal property \ref{prop.universal} and define an automorphism $\al_s$ of $\Ptilde$ satisfying $\al_s(x) = x$ for all $x \in \Mtilde$ and $\al_s(u) = u v_s$. Assume that $x_i \in \Mtilde$ and $\eps_1,\ldots,\eps_n \in \{1,-1\}$ such that $x_0 u^{\eps_1} x_1 \cdots u^{\eps_n} x_n$ is a reduced element in $\Ptilde$. We have to prove that
$$\tau\bigl(x_0 (u v_s)^{\eps_1} x_1 \cdots (u v_s)^{\eps_n} x_n \bigr) = 0 \; . $$
But grouping the $x_i$ with the $(v_s)^{\eps_j}$, one can consider $x_0 (u v_s)^{\eps_1} x_1 \cdots (u v_s)^{\eps_n} x_n$ as a reduced element in $\Ptilde$ as well. Therefore, its trace is zero.

{\bf Claim 1.} We have $\Ptilde = P *_N (N \ot \rL \Z)$. To prove this claim, we need to check that
$$\tau(x_0 v^{k_1} x_1 \cdots v^{k_n} x_n) = 0$$
whenever $x_0,x_n \in P$, $x_1,\ldots,x_{n-1} \in P \ominus N$ and $k_1,\ldots,k_n \in \Z - \{0\}$. To show this, we may assume that either $x_i \in M$ or that $x_i = y_{i,0} u^{\eps_{i,1}} y_{i,1} \cdots u^{\eps_i,\ell_i} y_{i,\ell_i}$ is a reduced element in $P$ of length $\ell_i \geq 1$. Make the following observations.
\begin{itemize}
\item If $k \in \Z - \{0\}$, then $M v^k M \subset \Mtilde \ominus M$.
\item If $1 \leq i \leq j \leq n-1$ and $x_i,x_{i+1},\ldots,x_j \in M \ominus N$, then
$$M v^{k_i} x_i v^{k_{i+1}} \cdots x_j v^{k_{j+1}} M \subset \Mtilde \ominus M \; .$$
\item We have that $\Mtilde \ominus M$ is a subspace of both $\Mtilde \ominus N$ and $\Mtilde \ominus \theta(N)$.
\end{itemize}
Altogether it follows that $x_0 v^{k_1} x_1 \cdots v^{k_n} x_n$ can be considered as a reduced element of $\HNN(\Mtilde,N,\theta)$ so that it indeed has trace zero.

\begin{remark}\label{rem.weakcont}
Note that claim 1 implies the following: if $N$ is amenable, then $\bim{P}{\rL^2(\Ptilde \ominus P)}{P}$ is weakly contained in the coarse $P$-$P$-bimodule. A detailed argument can be found e.g.\ in \cite[Proposition 3.1]{CH08}.
\end{remark}

Using claim 1 one checks easily that
$$E_P(\al_s(x)) = \rho_s^n \, x \quad\text{whenever}\;\; x = x_0 u^{\eps_1} x_1 \cdots u^{\eps_n} x_n \;\;\text{is a reduced element in $P$}\; .$$
It follows in particular that for all $0 < \rho < 1$, formula \eqref{eq.defmrho} yields a well defined normal unital completely positive map on $P$.

Also, since the automorphism $\beta \in \Aut \Mtilde$ is the identity on $M$, and so in particular on $N$ and $\theta(N)$, we can extend $\beta$ to an automorphism of $\Ptilde$ satisfying $\beta(u) = u$. By construction $\beta \circ \al_s = \al_{-s} \circ \be$ and $\beta(x) = x$ for all $x \in P$.

{\bf Claim 2.} We have $\Ptilde = P *_M \al_1(P)$. To prove this claim, first note that the von Neumann subalgebra of $\Ptilde$ generated by $P$ and $\al_1(P)$ contains $M$, $u$ and $\al_1(u) = uv$. Hence, it contains $M$, $u$ and $v$, which means that it equals $\Ptilde$. To conclude we have to prove that
$$\tau(x_0 \al_1(y_1) x_1 \cdots \al_1(y_n) x_n) = 0$$
whenever $x_0,x_n \in P$ and $x_1,\ldots,x_{n-1}, y_1,\ldots,y_n \in P \ominus M$. It is sufficient to prove this statement when all $x_i, y_j$ are reduced elements of $P$ with $x_1,\ldots,x_{n-1},y_1,\ldots,y_n$ having length at least $1$. Using the fact that $MvM$ and $M v^*M$ are subspaces of $\Mtilde \ominus M$, a straightforward computation gives the desired result.

With the machinery developed so far, we can now prove the following theorem. It is the HNN extension counterpart of \cite[Theorem 4.3]{IPP05} (see also \cite[Theorem 5.4]{PV09}). We will make use of the identifications in claims 1 and 2, resulting in a relatively elementary proof.

\begin{theorem}\label{thm.embedding}
Let $P = \HNN(M,N,\theta)$ be an HNN extension of finite von Neumann algebras, $p \in P$ a non-zero projection and $Q \subset pPp$ a von Neumann subalgebra. Denote by $S \subset pPp$ the quasi-normalizer of $Q$ inside $pPp$ (see Section \ref{subsec.amal} for terminology). Let $\m_\rho$ be the completely positive maps defined in \eqref{eq.defmrho}.

If there exist $0 < \rho < 1$ and $\delta > 0$ such that $\tau(d^* \m_\rho(d)) \geq \delta$ for all $d \in \cU(Q)$, then either $Q \embed_P N$ or $S \embed_P M$.
\end{theorem}
\begin{proof}
We assume that $Q \not\embed_P N$ and prove that $S \embed_P M$. Since $\tau(d^* \m_\rho(d))$ increases when $\rho$ increases towards $1$, we can take $s$ of the form $s=(2n)^{-1}$ and $\delta > 0$ such that
$$\tau(d^* \al_s(d)) = \tau\bigl( d^* E_P(\al_s(d))\bigr) = \tau(d^* \m_{\rho_s}(d)) \geq \delta \quad\text{for all}\;\; d \in \cU(Q) \; .$$
Denote by $y \in p \Ptilde \al_s(p)$ the unique element of minimal $\|\, \cdot\,\|_2$ in the weakly closed convex hull of $\{d^* \al_s(d) \mid d \in \cU(Q)\}$. It follows that $\tau(y) \geq \delta$ so that $y \neq 0$. By uniqueness of $y$, we have $d y = y \al_s(d)$ for all $d \in Q$. Denote by $w_0 \in p \Ptilde \al_s(p)$ the polar part of $y$. Then $w_0$ is a partial isometry satisfying $d w_0 = w_0 \al_s(d)$ for all $d \in Q$. Denote $p_0 := w_0 w_0^*$ and note that $p_0 \in p\Ptilde p \cap Q'$. By claim 1 above, we have $\Ptilde = P *_N (N \ot \rL \Z)$. Combined with Lemma \ref{lemma.IPP} and the assumption that $Q \not\embed_P N$, it follows that $p_0 \in p P p$. Similarly $w_0^* w_0 = \al_s(p_1)$ for some projection $p_1 \in p P p$. Since $\beta(x) = x$ for all $x \in P$ and $\be \circ \al_s = \al_{-s} \circ \be$, one checks that $w_1 := w_0 \al_{2s}(\be(w_0^*))$ is a partial isometry in $\Ptilde$ with left support projection $p_0$, right support projection $\al_{2s}(p_0)$ and $d w_1 = w_1 \al_{2s}(d)$ for all $d \in Q$. Putting
$$w = w_1 \, \al_{2s}(w_1) \, \al_{4s}(w_1) \, \cdots \, \al_{2(n-1)s}(w_1) \; ,$$
we have found a partial isometry $w \in \Ptilde$ with left support projection $p_0$, right support projection $\al_1(p_0)$ and $d w = w \al_1(d)$ for all $d \in Q$.

Since $p_0 \in pPp \cap Q'$, we have in particular that $p_0 \in S$. Again using claim 1 above, Lemma \ref{lemma.IPP} and the assumption that $Q \not\embed_P N$, the relation $d w = w \al_1(d)$ implies that $w^* S w = \al_1(p_0 S p_0)$. In particular $S \embed_{\Ptilde} \al_1(P)$. By claim 2 above we can view $\Ptilde$ as $\Ptilde = P *_M \al_1(P)$. By the remark following Lemma \ref{lemma.IPP} it then follows that $S \embed_P M$.
\end{proof}

\begin{remark}\label{rem.lengthdef}
By \cite{Ue07} amalgamated free products and HNN extensions are related to each other up to amplifications. We comment on this in Remark \ref{rem.ueda} and show there how our proof for Theorem \ref{thm.embedding}, which is rather simple thanks to claim 2, can be used to recover as well the results in \cite[Theorem 4.3]{IPP05} and \cite[Theorem 5.4]{PV09}.
\end{remark}

%%%%%%%%%%%%%%%%%%%%%%%%%%%%%%%%%%%%%%%%%%%%%%%%%%%%%%%%%%%%%%%%%%%%%%%%%%%%%%%%%
\section{Uniqueness of group measure space Cartan subalgebras}
%%%%%%%%%%%%%%%%%%%%%%%%%%%%%%%%%%%%%%%%%%%%%%%%%%%%%%%%%%%%%%%%%%%%%%%%%%%%%%%%%

We prove the following version of Theorem \ref{UniqueCartanHNN} stated in the introduction, also allowing for arbitrary amplifications.

\begin{theorem}\label{UniqueCartanHNN-amplif}
Let $\Gamma = \HNN(H,\Sigma,\theta)$ be an HNN extension satisfying all the hypotheses in Theorem \ref{UniqueCartanHNN}. Let $\Gamma \actson (X,\mu)$ be an arbitrary free ergodic p.m.p.\ action and denote $P := \rL^\infty(X) \rtimes \Gamma$.

If $\Lambda\curvearrowright (Y,\eta)$ is an arbitrary free ergodic p.m.p.\ action, $p \in \M_n(\C)\ot P$ is a projection and
$$\pi\,:\,\rL^{\infty}(Y)\rtimes\Lambda\rightarrow p(\M_n(\C)\ot P)p$$
is a $*$-isomorphism, then there exist a projection $q \in \D_n(\C)\ot \rL^{\infty}(X)$ and partial isometry $u\in \M_n(\C)\ot P$ such that $uu^* = q$, $u^* u=p$ and $$\pi(\rL^{\infty}(Y))=u^*(\D_n(\C)\ot \rL^{\infty}(X))u \; .$$
Here $\D_n(\C) \subset \M_n(\C)$ denotes the subalgebra of diagonal matrices.
\end{theorem}

To prove Theorem \ref{UniqueCartanHNN-amplif}, we first apply the following lemma to the action $\Gamma \actson \M_n(\C) \ot \rL^\infty(X)$. It follows that $\pi(\rL^\infty(Y)) \embed A \rtimes \Sigma$. Next we use the finiteness $\bigcap_{i=1}^n g_i \Sigma g_i^{-1}$ together with \cite[Theorem 6.16]{PV06}, to conclude that $\pi(\rL^\infty(Y)) \embed A$. Finally the conclusion of Theorem \ref{UniqueCartanHNN-amplif}  follows from \cite[Theorem A.1]{Po01}.

\begin{lemma}\label{lemma.crucial}
Let $H$ be a group that contains a non-amenable subgroup with the relative property (T) or that contains two commuting non-amenable subgroups. Let $\Sigma < H$ be an amenable subgroup and $\theta : \Sigma  \recht H$ an injective group homomorphism. Denote by $\Gamma = \HNN(H,\Sigma,\theta)$ the corresponding HNN extension. Let $\Gamma \actson (A,\tau)$ be any trace preserving action on the amenable von Neumann algebra $A$.

If $p \in A \rtimes \Gamma$ is a projection and $B \rtimes \Lambda = p(A \rtimes \Gamma)p$ is any other crossed product decomposition with $B$ abelian, then $B \embed A \rtimes \Sigma$.
\end{lemma}
\begin{proof}
Write $P = A \rtimes \Gamma$ and put $M = A \rtimes H$, $N = A \rtimes \Sigma$. Let $t \in \Gamma$ be the stable letter, with corresponding unitary $u_t \in \rL \Gamma \subset A \rtimes \Gamma$. Define the $*$-homomorphism $\thetatilde : N \recht M : \thetatilde(x) = u_t x u_t^*$ for all $x \in N$. Note that $\thetatilde(a u_g) = \sigma_t(a) u_{\theta(g)}$ for all $a \in A$ and $g \in \Sigma$. The universal property \ref{prop.universal} easily yields that
$P = \HNN(M,N,\thetatilde)$. So we have on $P$ the completely positive maps $\m_\rho$ and their dilation $\al_t$, as constructed in Section \ref{sec.deform}.

We apply the transfer of rigidity Lemmas 3.1 and 3.2 in \cite{PV09} to the amenable subalgebra $N \subset P$. Note that $P$ indeed contains either a non-amenable von Neumann subalgebra with the relative property (T) or two commuting non-amenable von Neumann subalgebras. Using the inequality $\|x-\m_{\rho_t}(x)\|_2 \leq \|\alpha_t(x)-E_P(\alpha_t(x))\|_2$ for all $x \in P$, we get from \cite[Lemmas 3.1 and 3.2]{PV09} a $0 < \rho < 1$ and a sequence $(s_k)_k\in\Lambda$ such that
\begin{itemize}
\item $\|\m_{\rho}(v_{s_k})-s_k\|_2\leq \tau(p)/5000$ for all $k$,
\item $\|E_{N}(xv_{s_k}y)\|_2\rightarrow 0$ for all $x,y\in P$.
\end{itemize}
Lemma \ref{lemma.end} below provides a $0<\rho<1$ and a $\delta>0$ such that $\tau(d^* \m_{\rho}(d))\geq \delta$ for all $d\in\mathcal{U}(B)$. By Theorem \ref{thm.embedding} we either reach the conclusion of the lemma, or we find that the normalizer of $B$ embeds into $M$ inside $P$. Since the normalizer of $B$ is the whole of $p P p$, the latter is absurd.
\end{proof}

We do not give a proof for the following result, since it is, mutatis mutandis, the same as the proof of \cite[Lemma 5.7]{PV09}. Two observations have to be made. First, as we recalled in Section \ref{HNNgrp}, every HNN extension $\Gam$ can act on its Bass-Serre tree such that the natural length function on $\Gam$ is given by the geodesic distance on the tree. As a consequence, \cite[Lemma 4.1]{PV09} remains valid for HNN extensions. Secondly, we need a replacement for the combinatorial lemma \cite[Lemma 5.5]{PV09}. We provide this replacement as Lemma \ref{combinatoire} below.

\begin{lemma}\label{lemma.end}
Let $\Gamma = \HNN(H,\Sigma,\theta)$ be any HNN extension and $\Gamma \actson (A,\tau)$ any trace preserving action. Consider the completely positive maps $\m_\rho$ on $A \rtimes \Gamma$ as above. Let $B \subset p(A \rtimes \Gamma)p$ be an abelian von Neumann subalgebra. Assume that there exists $0 < \rho < 1$ and a sequence of unitaries $v_k \in p(A \rtimes \Gamma)p$ that normalize $B$ and satisfy
\begin{itemize}
\item $\|v_k - \m_\rho(v_k)\|_2 \leq \tau(p)/5000$ for all $k$,
\item $\|E_{A \rtimes \Sigma}(x v_k y)\|_2 \recht 0$ for all $x,y \in A \rtimes \Gamma$.
\end{itemize}
Then there exists $0 < \rho_0 < 1$ and a $\delta > 0$ such that $\tau(d^* \m_{\rho_0}(d)) \geq \delta$ for all $d \in \cU(B)$.
\end{lemma}

To conclude we provide the promised combinatorial lemma. Let $\Gamma = \HNN(H,\Sigma,\theta)$ be an HNN extension and $\Gamma \actson (A,\tau)$ be a trace preserving action. Use the notations introduced in Section \ref{HNNgrp}. Let $g\in\Gam$ be an element of length at least one and write a reduced expression $g=x_0t^{\eps_1}\cdots t^{\eps_n}x_n$ with $x_i\in H$ and $\eps_i\in\{-1,1\}$. We call $x_0t^{\eps_1}$ the \textit{first letter} of $g$ and $t^{\eps_n}x_n$ the \textit{last letter} of $g$.

Let $\eps,\eps'\in\{-1,1\}$ and $g'\in Ht^{\eps}$, $h'\in t^{\eps'}H$. Let $\Gam_{g' ,h'}\subset\Gam$ be the set of elements in $\Gam$ beginning by $g'$ and ending by $h'$. This means that $g\in\Gam_{g' ,h'}$ if and only if every reduced expression for $g$ begins with $g'\sigma$ with some $\sigma\in\Sigma_{\eps}$ and ends with $\sigma'h'$ with some $\sigma'\in\Sigma_{-\eps'}$.

We denote by $\mathcal{P}_{g',h'}$ the orthogonal projection of $\rL^2(P)$ onto the closed linear span of $\{au_g\,|\,g\in\Gam_{g',h'},a\in A\}$.
In general, whenever $W\subset\Gam$, we denote by $\mathcal{P}_W$ the orthogonal projection of $\rL^2(P)$ onto the closed linear span of $\{au_g\,|\,g\in W,a\in A\}$.
We also denote by $\mathcal{P}_K$ the orthogonal projection of $\rL^2(P)$ onto the closed linear span of $\{au_g\,|\,|g|\leq K,a\in A\}$.

\begin{lemma}\label{combinatoire}
Let $K\in\N-\{0,1\}$ and assume that $(y_k)$ is a bounded sequence in $P$ with the following properties.
\begin{itemize}
\item $y_k=\mathcal{P}_K(y_k)$ for all $k$,
\item $\|E_{N}(xy_kz)\|_2\rightarrow 0$ for all $x,z\in P$.
\end{itemize}
Let $g,h\in\Gam$ with $|g|,|h|\geq K+1$ and write $g,h$ as reduced expressions. Denote by $g'$ the first letter of $g$ and by $h'$ the last letter of $h$. Then, we can write
$$u_{g}y_ku_{h}=a_k+b_k$$
where $a_k,b_k$ are bounded sequences in $P$ satisfying the following properties.
\begin{itemize}
\item $a_k=\mathcal{P}_{g',h'}(a_k)$ for all $k$,
\item $\|b_k\|_2\rightarrow 0$.
\end{itemize}
\end{lemma}
\begin{proof}
Choose reduced expressions for $g$ and $h$:
$$ g=x_0t^{\eps_1}\cdots t^{\eps_n}x_n\quad\text{and}\quad h=y_0t^{\eta_1}\cdots t^{\eta_m}y_m,$$
where $x_i,y_i\in H$, $\eps_i,\eta_i\in\{-1,1\}$ and $n,m\geq K+1$. For $0\leq i\leq n$ and $0\leq j\leq m$ define
$$g_i=x_{n-i}t^{\eps_{n-i+1}} \cdots x_{n-1} t^{\eps_n} \; x_n\quad\text{and}\quad h_j=y_0 \; t^{\eta_1}y_1 \cdots t^{\eta_j}y_j \; .$$
We use the convention $g_0=e=h_0$. Observe that $|g_i|=i$ and $|h_j|=j$.

Write $N_{\eps}=A\rtimes\Sigma_{\eps}$ for $\eps\in\{-1,1\}$.
For  $0\leq i\leq n$ and $0\leq j\leq m$ we define $W_{i,j}=gg_i^{-1}\Sigma_{\eps_{n-i}}h_j^{-1}h$. Put $W=\cup_{i+j\leq K}W_{i,j}$ and observe that
$$\mathcal{P}_{W_{i,j}}(a)=u_{gg_i^{-1}}E_{N_{\eps_{n-i}}}(u_{g_ig^{-1}} a u_{hh_j^{-1}})u_{h_j^{-1}h}\quad\text{for all}\,\,\,a\in P.$$
Hence, $\mathcal{P}_{W_{i,j}}$ is bounded as a map from $P$ to $P$. The orthogonal projections $\mathcal{P}_{W_{i,j}}$ commute and hence
$$1-\mathcal{P}_{W}=\prod_{i+j\leq K}(1-\mathcal{P}_{W_{i,j}}).$$
is bounded on $P$. We put $b_k=\mathcal{P}_{W}(u_{g}y_ku_{h})$, $a_k=u_{g}y_ku_{h}-b_k$. So, $a_k$ and $b_k$ are bounded sequences in $P$ with $u_{g}y_ku_{h}=a_k+b_k$.

Since $\theta(N) = uNu^*$, it follows that $\|E_{N_\eps}(xy_kz)\|_2\rightarrow 0$ for all $x,z\in P$ and all $\eps\in\{-1,1\}$. It follows that
$$\|b_k\|_2^2\leq\sum_{i+j\leq K}\|\mathcal{P}_{W_{i,j}}(u_{g}y_ku_{h})\|_2^2=\sum_{i+j\leq K}\|E_{N_{\eps_{n-i}}}(u_{g_i}y_ku_{h_j^{-1}})\|_2^2\rightarrow 0.$$
It remains to prove that $a_k=\mathcal{P}_{g',h'}(a_k)$ for all $k$. But, if $r\in\Gam$ with $|r|\leq K$ and if $grh$ admits no reduced expression that begins with $g'$ and ends with $h'$, there must exist $i,j$ with $i+j\leq K$ such that $g_irh_j\in\Sigma_{\eps_{n-i}}$. Hence, $grh\in W$. As a consequence, whenever $y\in M$ with $y=\mathcal{P}_{K}(y)$, we have
$$u_{g}yu_{h}-\mathcal{P}_W(u_{g}yu_{h})\in\mathcal{P}_{g',h'}(\rL^2(P)).$$
This concludes the proof of the lemma.
\end{proof}

\begin{remark}\label{rem.rigidinGamma}
In the formulation of Theorem \ref{UniqueCartanHNN} and its amplified version \ref{UniqueCartanHNN-amplif}, we could as well assume that instead of $H$, the larger group $\Gamma = \HNN(H,\Sigma,\theta)$ contains a non-amenable subgroup with the relative property (T) or contains two commuting non-amenable subgroups. The proof of the theorem remains unchanged. On the other hand, this is not a real generalization, because such rigid subgroups have to lie inside $H$ (after conjugacy and passage to a finite index subgroup). Nevertheless this remark will be useful in the proof of Theorem \ref{UniqueCartanTrees}. A similar remark applies to \cite[Theorem 1.1]{PV09}.
\end{remark}

\begin{remark}\label{rem.ueda}
We are grateful to the referee of the first version of this article who pointed us towards the following result in \cite{Ue07}, inspired by similar observations for equivalence relations in \cite{Ga99,Pa99}. Let $M$ be a tracial von Neumann algebra, $N \subset M$ a von Neumann subalgebra and $\theta : N \recht M$ a trace preserving embedding. Consider the trace preserving embeddings
$$N \oplus N \hookrightarrow \M_2(\C) \ot M : x \oplus y \mapsto \begin{pmatrix} x & 0 \\ 0 & \theta(y) \end{pmatrix}\quad\text{and}\quad N \oplus N \hookrightarrow \M_2(\C) \ot N : x \oplus y
\mapsto \begin{pmatrix} x & 0 \\ 0 & y \end{pmatrix} \; .$$
Let $u \in \HNN(M,N,\theta)$ be the stable unitary and denote by $(e_{ij})$, resp.\ $(f_{ij})$ the canonical matrix units in $\M_2(\C) \ot M$, resp.\ $\M_2(\C) \ot N$.
By \cite[Proposition 3.1]{Ue07} there is a canonical trace preserving $*$-isomorphism
$$\Psi : \HNN(M,N,\theta) \recht e_{11} \bigl((\M_2(\C) \ot M) \underset{N \oplus N}{*} (\M_2(\C) \ot N)\bigr)e_{11} : \begin{cases} \Psi(x) = e_{11} x \;\;\text{for all}\;\; x \in M \;, \\ \Psi(u) = e_{12} f_{21} \; ,\end{cases}$$
where the amalgamated free product is with respect to the embeddings above and the unique trace preserving conditional expectations.

Since on the group level HNN extensions cannot be canonically written as amalgamated free products, one cannot directly deduce Theorem \ref{UniqueCartanHNN-amplif} from the analogous \cite[Theorem 5.2]{PV09}. Nevertheless the $*$-isomorphism $\Psi$ can be used as the starting point for an alternative proof for Theorem \ref{UniqueCartanHNN-amplif} that we sketch now.
One actually has to generalize \cite[Theorem 5.6]{PV09} from crossed products with amalgamated free product groups to the following more general statement about arbitrary amalgamated free products.

{\it Suppose that $M = M_1 *_P M_2$ is an amalgamated free product of tracial von Neumann algebras w.r.t.\ the unique trace preserving conditional expectations. Assume that $P$ is amenable and that $M$ admits a von Neumann subalgebra $M_0$ without amenable direct summand such that either $M_0 \subset M$ has the relative property (T) or $M_0' \cap M$ has no amenable direct summand. If $M = B \rtimes \Lambda$ for some abelian von Neumann algebra $B \subset M$ and some countable group $\Lambda$ acting on $B$, then $B \embed_M P$.}

Such a generalization is not totally innocent and more technical than the proof that we gave above for Theorem \ref{UniqueCartanHNN-amplif}. One first has to replace the usage of Herz-Schur multipliers in \cite[Section 4]{PV09} by the corresponding results in \cite[Theorem 3.1 and Section 5]{RX05}. Next one should redo the proofs of \cite[Lemmas 5.5 and 5.7]{PV09} within the general framework of amalgamated free products of von Neumann algebras.

Conversely, by \cite[Proposition 3.4]{Ue07} also amalgamated free products can essentially be viewed as an HNN extension. So, let $M_1$ and $M_2$ be tracial von Neumann algebras and let $P$ be a common von Neumann subalgebra on which the traces coincide. We consider the amalgamated free product $M_1 *_P M_2$ w.r.t.\ the trace preserving conditional expectations. Denote $M = M_1 \oplus M_2$ and define $N \subset M$ given by $N = P \oplus P$. Put $\theta : N \recht M : \theta(x \oplus y) = y \oplus x$. Let $\HNN(M,N,\theta)$ be the HNN extension with stable unitary $u$ and put $p = 1 \oplus 0$ One checks easily that there is a unique trace preserving $*$-homomorphism
$$\Theta : M_1 *_P M_2 \recht p \HNN(M,N,\theta) p : \Theta(x) = \begin{cases} x \oplus 0 &\;\;\text{if}\;\; x \in M_1 \; , \\ u (0 \oplus x) u^* &\;\;\text{if}\;\; x \in M_2 \; .\end{cases}$$
More precisely, $\Theta$ extends to a $*$-isomorphism between $M_1 *_P M_2 *_P (P \ot \rL \Z)$ and $p \HNN(M,N,\theta) p$ sending the extra generator $v \in \rL \Z$ to $u^2 p$.

The $*$-isomorphisms $\Psi$ and $\Theta$ intertwine the length deformation on the amalgamated free product and the HNN extension (see Section \ref{sec.deform}). So, Theorem \ref{thm.embedding} and \cite[Theorem 5.4]{PV09} can be deduced from each other. As explained in Remark \ref{rem.lengthdef}, our proof of Theorem \ref{thm.embedding} is slightly simpler which is our reason to present it in full detail in Section \ref{sec.deform}.
\end{remark}

%%%%%%%%%%%%%%%%%%%%%%%%%%%%%%%%%%%%%%%%%
\section{Proof of Theorem \ref{UniqueCartanTrees}}
%%%%%%%%%%%%%%%%%%%%%%%%%%%%%%%%%%%%%%%%%

In Theorem \ref{UniqueCartanHNN} we have seen that the conclusion of Theorem \ref{UniqueCartanTrees} holds for certain HNN extensions. In \cite[Theorem 1.1]{PV09} it was shown that the conclusion also holds for certain amalgamated free product groups. So, it suffices to prove that all groups satisfying the assumptions of Theorem \ref{UniqueCartanTrees} fall into one of both families.

Take $\Gamma \actson \cT$ with the properties assumed in Theorem \ref{UniqueCartanTrees}. Let $e \in \rE(\cT)$ be an edge such that the stabilizer $\Sigma := \Stab e$ is amenable and such that the smallest subtrees containing $\Gamma \cdot s(e)$, resp.\ $\Gamma \cdot r(e)$, are both equal to the whole of $\cT$. We claim that there exist $g_1,\ldots,g_m \in \Gamma$ such that $\bigcap_{i=1}^m g_i \Sigma g_i^{-1}$ is finite. Let $\cT_0 \subset \cT$ be a finite subtree with finite stabilizer. Denote by $p_1,\ldots,p_n$ the vertices of $\cT_0$. Since the smallest subtree containing all the edges in $\Gamma \cdot e$ is the whole of $\cT$, we can take $h_i,k_i \in \Gamma$ such that $p_i$ lies on the geodesic path joining $g_i \cdot e$ to $h_i \cdot e$. It is easy to check that
$$\bigcap_{i=1}^n \bigl( g_i \Sigma g_i^{-1} \cap h_i \Sigma h_i^{-1} \bigr) \subset \Stab \cT_0 \; ,$$
proving the claim.

Using Remark \ref{rem.rigidinGamma} it remains to prove that $\Gamma$ is either a non-trivial amalgamated free product over $\Sigma$ or an HNN extension over $\Sigma$. By \cite{Se83} we know that the quotient graph $\cG := \cT / \Gamma$ can be equipped with the structure of a graph of groups $(\cG,\{\Gamma_q\}_{q \in \rV(\cG)}, \{\Sigma_e\}_{e \in \rE(\cG)})$ such that $\Gamma$ is the fundamental group of this graph of groups and $\cT$ is its Bass Serre tree. For every $e \in \rE(\cG)$, denote by $s_e : \Sigma_e \recht \Gamma_{s(e)}$ and $r_e : \Sigma_e \recht \Gamma_{r(e)}$ the structural injective group homomorphisms. Denote by $\pi : \cT \recht \cG = \cT / \Gamma$ the quotient map.
We fixed in the previous paragraph a favorite edge $e$ and also denote by $e$ its image in $\cG = \cT/\Gamma$. We may assume that $\Sigma_e = \Sigma$. Note also that to every connected subgraph $\cG_1 \subset \cG$ corresponds a subgroup $\Gamma_1 < \Gamma$, given as the fundamental group of our graph of groups restricted to $\cG_1$, as well as a $\Gamma_1$-invariant subtree $\cT_1 \subset \cT$ satisfying $\pi(\cT_1) = \cG_1$.

Consider the graph $\Gr'$ obtained from the graph $\Gr$ by removing the edges $e$ and $\bar{e}$. There are two cases.

\textbf{Case 1.} The graph $\Gr'$ is connected. We will contract the connected subgraph $\Gr'\subset\Gr$ in one vertex and use Serre's ``d\'evissage'' technique to conclude. Following \cite{Se83}, we define $H$ to be the fundamental group of our graph of groups restricted to $\Gr'$. Via the source homomorphism $s_e$ we view $\Sigma = \Sigma_e$ as a subgroup of $H$. Define the injective group homomorphism $\theta : \Sigma \recht H$ given by $r_e$. It follows from \cite[Lemme 6, Section 5.2]{Se83} that $\Gam\cong\HNN(H,\Sigma,\theta)$.

\textbf{Case 2.} The graph $\Gr'$ is not connected. Let $\Gr_1$ (resp. $\Gr_2$) be the connected component of $s(e)$ (resp.\ $r(e)$). We use a two step contraction procedure. Define for $i=1,2$, the subgroup $\Gam_i < \Gamma$ as the fundamental group of our graph of groups restricted to the connected subgraph $\Gr_i$ of $\Gr$.
By the source homomorphism $s_e$ we view $\Sigma = \Sigma_e$ as a subgroup of $\Gamma_1$ and by the range homomorphism $r_e$ we view $\Sigma$ as a subgroup of $\Gamma_2$. Using twice \cite[Lemme 6, Section 5.2]{Se83}, first contracting $\Gr_1$ to one vertex and then contracting $\Gr_2$ to one vertex, we conclude that $\Gam\cong \Gam_1 *_\Sigma \Gam_2$. Finally observe that we decomposed $\Gamma$ as a non-trivial amalgamated free product: $\Gamma_1 \neq \Sigma \neq \Gamma_2$. Indeed, otherwise $\Gamma = \Gamma_1$ or $\Gamma = \Gamma_2$. If $\Gamma = \Gamma_1$, the subtree $\cT_1 \subset \cT$ corresponding to the connected subgraph $\cG_1 \subset \cG$, is $\Gamma$-invariant and contains $s(e)$. By our assumptions, it follows that $\cT_1 = \cT$. Since $\pi(\cT_1) = \cG_1$, this is a contradiction. The equality $\Gamma = \Gamma_2$ leads to a contradiction in a similar way.

%%%%%%%%%%%%%%%%%%%%%%%%%%%%%%%%%%%%%%%%%%%%%%%%%%%%%%%%%%%%%%%%%%%%%%%%%%%%%%%%%%%%%%%%%%
\section{Stable W$^*$-superrigidity}
%%%%%%%%%%%%%%%%%%%%%%%%%%%%%%%%%%%%%%%%%%%%%%%%%%%%%%%%%%%%%%%%%%%%%%%%%%%%%%%%%%%%%%%%%%

Recall that a free ergodic p.m.p.\ action $\Gamma \actson (X,\mu)$ is said to be W$^*$-superrigid if the following holds: if $\Lambda \actson (Y,\eta)$ is any other free ergodic p.m.p.\ action and $\theta : \rL^\infty(Y) \rtimes \Lambda \recht \rL^\infty(X) \rtimes \Gamma$ is a $*$-isomorphism, then the actions $\Lambda \actson (Y,\eta)$ and $\Gamma \actson (X,\mu)$ are conjugate and the isomorphism $\theta$ is in a precise sense implemented by this conjugacy and a scalar $1$-cocycle (see \cite[Definition 6.1]{PV09}).

A somehow more natural notion than W$^*$-superrigidity is stable W$^*$-superrigidity, appropriately taking care of amplifications. This is discussed in detail in \cite[Section 6.2]{PV09}. We only mention here that stable W$^*$-superrigidity for $\Gamma \actson (X,\mu)$ implies W$^*$-superrigidity when $\Gamma$ has no non-trivial finite normal subgroups and if finite index subgroups of $\Gamma$ act ergodically on $(X,\mu)$.

Theorem \ref{thm.Wstarsuperrigid} is a particular case of the following result.

\begin{theorem}\label{thm.Wstarsuperrigid-general}
Let $H$ be a countable group with infinite amenable subgroup $\Sigma < H$. Let $\theta : \Sigma \recht H$ be an injective group homomorphism. Denote by $\Gam=\HNN (H,\Sigma,\theta)$ the HNN extension and suppose that there exist $g_1,\ldots, g_n\in \Gam$ such that $\bigcap_{i=1}^n g_i\Sigma g_i^{-1}$ is finite. Note that this last condition is automatic if $\Sigma\cap\theta(\Sigma)$ is finite.
\begin{enumerate}
\item If $H$ admits a non-amenable normal subgroup $H_0$ with the relative property $(T)$, then all of the following actions are stably W$^*$-superrigid.
\begin{itemize}
\item Every free p.m.p.\ action $\Gam\curvearrowright (X,\mu)$ whose restriction to $H$ is a generalized Bernoulli action $H\curvearrowright (X_0,\mu_0)^I$ with the property that $H_0 \cdot i$ and $\Sigma \cdot i$ are infinite for all $i\in I$.
\item Every free p.m.p.\ action $\Gam\curvearrowright (X,\mu)$ whose restriction to $H$ is a Gaussian action defined by an orthogonal representation $\pi\, :\, H\rightarrow \mathcal{O}(\mathcal{K}_{\R})$ with the property that the restrictions $\pi_{| H_0}$ and $\pi_{| \Sigma}$ have no non-zero finite dimensional subrepresentations.
\end{itemize}
\item Suppose that $H$ admits non-amenable commuting subgroups $H_0$ and $H_1$ such that $H_0$ is normal in $H$. If $\Gam\curvearrowright (X,\mu)$ is a free p.m.p.\ action whose restriction to $H$ is a generalized Bernoulli action $H\curvearrowright (X_0,\mu_0)^I$ with the property that $H_1 \cap \Stab i$ is amenable for all $i\in I$ and that $H_0 \cdot i$ and $\Sigma \cdot i$ are infinite for all $i\in I$, then $\Gam\curvearrowright (X,\mu)$ is stably W$^*$-superrigid.
\end{enumerate}
\end{theorem}

\begin{proof}
By the uniqueness of group measure space Cartan Theorem \ref{UniqueCartanHNN-amplif} and using \cite[Lemma 6.5]{PV09}, it is sufficient to prove that for all group actions $\Gamma \actson (X,\mu)$ appearing in the theorem, we have that all measurable $1$-cocycles $\om : \Gamma \times X \recht \cG$ with values in either a countable group $\cG$ or the group $\cG = S^1$, are cohomologous to a group morphism from $\Gamma$ to $\cG$. Take such a $1$-cocycle $\om$. In case 1 we apply \cite[Theorem 0.1]{Po05} and in case 2 we apply \cite[Theorem 1.1]{Po06a}. In both cases it follows that the restricted $1$-cocycle $\om_{|H \times X}$ is cohomologous to a group morphism from $H$ to $\cG$. We may therefore assume that $\om(g,x) = \delta(g)$ for all $g \in H$ and a.e.\ $x \in X$. Denote by $t \in \HNN(H,\Sigma,\theta)$ the stable letter. Because $t^{-1}Ht\cap H=\Sigma$ and $\Sigma\curvearrowright (X,\mu)$ is weakly mixing, \cite[Proposition 3.6]{Po05} implies that $\om$ is also independent of $x$ in $t$. Hence, $\om(g,x)$ is independent of $x$ for all $g \in \Gamma$.
\end{proof}

%%%%%%%%%%%%%%%%%%%%%%%%%%%%%%%%%%%%%%%%%%%%%%%%%%%%%%%%%%%%%%%%%%%%%%%%%%%%%%%%%%%%%%

\end{document}